\documentclass[10pt,journal,compsoc]{IEEEtran}
\usepackage{graphicx}
\usepackage{amsmath}
\usepackage{amssymb}
\usepackage{color}
\ifCLASSOPTIONcompsoc
  \usepackage[nocompress]{cite}
\else
  \usepackage{cite}
\fi
\ifCLASSINFOpdf
\else
\fi
\begin{document}
%
\title{Hitting times and resistance distances  of $q$-triangulation graphs: Accurate results and applications}

\author{Yibo~Zeng and~Zhongzhi~Zhang
\IEEEcompsocitemizethanks{\IEEEcompsocthanksitem Yibo Zeng and Zhongzhi Zhang are with the Shanghai Key Laboratory of Intelligent Information
Processing, School of Computer Science, Fudan University, Shanghai 200433, China.\protect\\
E-mail: zhangzz@fudan.edu.cn}
\thanks{}
}

%
%

\markboth{IEEE TRANSACTIONS ON NETWORK SCIENCE AND ENGINEERING, VOL. 0, NO. 0, MONTH YYYY}%
{Shell \MakeLowercase{\textit{et al.}}: Bare Demo of IEEEtran.cls for Computer Society Journals}
%




\IEEEtitleabstractindextext{%
\begin{abstract}
Graph operations or products, such as triangulation and Kronecker product have been extensively applied to model complex networks with striking properties observed in real-world complex systems. In this paper, we study hitting times and resistance distances of $q$-triangulation graphs. For a simple connected graph $G$, its $q$-triangulation graph $R_q(G)$ is obtained from $G$ by performing the $q$-triangulation operation on $G$. That is, for every edge $uv$ in $G$, we add $q$ disjoint paths of length $2$, each having $u$ and $v$ as its ends. We first derive the eigenvalues and eigenvectors of normalized adjacency matrix of $R_q(G)$, expressing them in terms of those associated with $G$. Based on these results, we further obtain some interesting quantities about random walks and resistance distances for $R_q(G)$, including two-node hitting time, Kemeny's constant, two-node resistance distance, Kirchhoff index, additive degree-Kirchhoff index, and multiplicative degree-Kirchhoff index. Finally, we provide exact formulas for the aforementioned quantities of iterated $q$-triangulation graphs, using which we provide closed-form expressions for those quantities corresponding to a class of  scale-free small-world graphs, which has been applied to mimic complex networks.
\end{abstract}

}

\maketitle

\IEEEdisplaynontitleabstractindextext

%
\IEEEpeerreviewmaketitle



%
%
%
%
%
%
%

\IEEEraisesectionheading{\section{Introduction}\label{sec:introduction}}
\IEEEPARstart{G}{raph} operations and products play an important role in network science, which have been used to model complex networks with the prominent scale-free~\cite{BaAl99} and small-world~\cite{WaSt98} properties as observed in various real-life networks~\cite{Ne03}. Since diverse realistic large-scale networks consist of smaller pieces or patterns, such as communities~\cite{GiNe02}, motifs~\cite{MiShItKaKhAl02}, and cliques~\cite{Ts15}, graph operations and products are a natural way to generate a massive graph out of smaller ones. Furthermore, there are many advantages to using graph operations and products to create complex networks. For example, it allows  analytical treatment for structural and dynamical aspects of the resulting networks. Thus far, a variety of graph operations and products have been introduced or proposed to construct models of complex networks, including triangulation~\cite{DoGoMe02,ZhRoZh07}, Kronecker product~\cite{We62,LeFa07,LeChKlFaGh10}, hierarchical product~\cite{BaCoDaFi09,BaDaFiMi09,BaCoDaFi16}, as well as corona product~\cite{LvYiZh15,ShAdMi17,QiLiZh18}.

Among  various graph operations and products, triangulation is a popular one. Let $G$ be a simple graph. The triangulation graph of $G$, denoted by $R(G)$, is the graph obtained from $G$ by performing triangulation operation~\cite{CvDoSa80}. That is, for each edge $uv$ in $G$, a new node $x$ is created and connected to nodes $u$ and $v$. Algebraic and combinatorial properties of triangulation graphs have been comprehensively studied~\cite{Ya06,XiZhCo16b,LiLiCa16,ShLiZh17,ShLiZh18}. For more convenient and practical applications, an extended triangulation operation called $q$-triangulation was proposed~\cite{ZhRoZh07, RoHaBe07}.  For a positive integer $q$, the $q$-triangulation graph of $G$, denoted by $R_q(G)$, is the graph obtained from $G$ by adding, for each edge $uv$ in $G$, $q$ disjoint paths of length $2$: $ux_1v$, $ux_2v$, $\ldots$, $ux_qv$. The $q$-triangulation operation can be iteratively used to a triangle, generating a model for complex networks with the scale-free small-world characteristics~\cite{ZhRoZh07, RoHaBe07}. However, the properties of $R_q(G)$ for a generic graph $G$ are not well understood, comparing to the traditional triangulation graph, i.e., $1$-triangulation graph.

In this paper, we provide an in-depth study on the properties for $q$-triangulation graph $R_q(G)$ of an arbitrary simple connected graph $G$, focusing on random walks and resistance distances, both of which have found wide applications~\cite{Gr06,DoSibu18}. We first give explicit formulas for eigenvalues and eigenvectors of normalized adjacency (or Laplacian) matrix for $R_q(G)$, based on which we determine two-node hitting time and the Kemeny's constant for random walks on $R_q(G)$. Also, we derive the expressions for two-node resistance distance, Kirchhoff index, additive degree-Kirchhoff index, and multiplicative degree-Kirchhoff index for $R_q(G)$. All obtained quantities for $R_q(G)$ are expressed in terms of those associated with $G$. Finally, we obtain explicit expressions for the aforementioned quantities of iterated $q$-triangulation of a graph $G$, and apply such results to a category of scale-free small-world networks~\cite{ZhZhZo07}, yielding closed-form formulas for several interesting quantities.

\section{Preliminaries}

In this section, we introduce some basic concepts for a graph, random walks and electrical networks.

\subsection{Graph and Matrix Notation}

Let $G(V,E)$ be a simple connected graph with $n$ nodes/vertices and $m$ edges.  The $n$ nodes form node  set $V(G)=\{1,2,\ldots,n\}$, and the $m$ edges constitute edge set $E(G)=\{e_1,e_2,\ldots, e_m\}$.

Let $A$ denote the adjacency matrix of $G$, whose entry $A(i,j)$ is 1 (or 0) if nodes $i$ and $j$ are (not) directly connected in $G$. Let $\Gamma (i)$ denote the set of neighbors of  node $i$  in graph $G$. Then the degree of node $i$ is  $d_i=\sum_{j  \in \Gamma (i)}A(i,j)$, which constitutes the $i$th entry of the diagonal degree matrix $D$ of $G$. The incidence matrix of  $G$ is an $n\times m$ matrix $B$, where $B(i,j)=1$ (or 0) if node $i$ is (not) incident with $e_j$. 

\newtheorem{lemma}{Lemma}
\begin{lemma}\label{bRank}~\cite{CvDoSa80}
Let $G$ be a simple connected graph with $n$ nodes. Then the rank of its incidence matrix $B$ is ${\rm rank} (B)=n-1$ if $G$ is bipartite, and ${\rm rank} (B)=n$  otherwise.
\end{lemma}

\begin{lemma}\label{BBAD}
Let $G$ be a simple connected graph. Then its incidence matrix $B$, adjacency matrix $A$ and diagonal degree matrix $D$ satisfy
\begin{equation*}
    BB^\top = A + D.
\end{equation*}
\end{lemma}

\subsection{Random Walks on Graphs}

For a graph $G$, we can define a discrete-time unbiased random walk taking place on it. At any time step, the walker jumps from its current location, node $i$,  to an adjacent node $j$ with probability $A(i,j)/d_i$. Such a random walk on $G$ is in fact a Markov chain~\cite{KeSe60} characterized by the transition probability matrix $T=D^{-1}A$, with its entry  $T(i,j)$ equal to $A(i,j)/d_i$. For a random walk on graph $G$, the stationary distribution is an $n$-dimension vector $\pi=(\pi_1, \pi_2, \ldots, \pi_n)$ satisfying $\pi  T=\pi$ and $\sum_{i=1}^n \pi_i=1 $. It is easy to verify that $\pi=(d_1/2m, d_2/2m, \ldots, d_n/2m)$ for unbiased random walks on graph $G$.
	
The transition  probability matrix $T$ of graph $G$ is generally asymmetric. However, $T$ is similar to the normalized adjacency matrix $P$ of $G$, which is defined by
\begin{equation*}
P=D^{-\frac{1}{2}} A D^{-\frac{1}{2}}=D ^{\frac{1}{2}}T D^{-\frac{1}{2}}\,.
\end{equation*}
Obviously, $P$ is symmetric, with the $(i,j)$th entry being $P(i,j)=\frac{A(i,j)}{\sqrt{d_id_j}}$. $I-P$ is the normalized Laplacian matrix~\cite{Ch97} of graph $G$, where $I$ is the $n\times n$ identity matrix.

\begin{lemma}\label{MinEig}\cite{Ch97}
Let $G$ be a simple connected graph with $n$ nodes, and let
$1=\lambda_1>\lambda_2\geq\ldots\geq\lambda_n\geq-1$ be  the eigenvalues of its normalized adjacency matrix $P$. Then
$\lambda_n=-1$ if and only if $G$ is bipartite.
\end{lemma}

Let $v_1$, $v_2$,$\ldots$, $v_n$ be the orthonormal eigenvectors corresponding to  the $n$ eigenvalues $\lambda_1$, $\lambda_2$,$\ldots$, $\lambda_n$, where $v_i=(v_{i1},v_{i2},\ldots,v_{in})^\top$. Then,
\begin{equation}\label{eig=1}
  v_1=\Big(\sqrt{d_1/2m},\sqrt{d_2/2m},...,\sqrt{d_n/2m}\Big)^\top
\end{equation}
and
\begin{equation}
\sum_{k=1}^n v_{ik}v_{jk}=\sum_{k=1}^n v_{ki}v_{kj}=\left\{
                                                      \begin{array}{ll}
                                                        1, & \hbox{if $i=j$;} \\
                                                        0, & \hbox{otherwise.}
                                                      \end{array}
                                                    \right.
\end{equation}
As for a bipartite graph $G$, whose node set $V(G)$ can be divided into two disjoint subsets $V_1$ and $V_2$, i.e., $V(G)=V_1\cup V_2$, we have
\begin{equation}\label{eig=-1}
  v_{ni}=\sqrt{d_i/2m},\, i\in V_1;~~v_{ni}=-\sqrt{d_i/2m}, \,  i\in V_2.
\end{equation}

A fundamental quantity related to random walks is  hitting time. The hitting time $T_{ij}$ from one node $i$ to another node $j$ is the expected number of jumps needed for a walker to first reach node $j$ starting from node $i$, which is relevant in various scenarios~\cite{Re01}. Many interesting quantities of graph $G$ can be defined or derived from hitting times.  For example, for a  graph $G$, its Kemeny's constant $K(G)$  is defined as  the expected number of steps required for a walker starting from node $i$ to a destination node, which is chosen randomly according to a stationary distribution of  random walks on $G$~\cite{Hu14}.  The Kemeny's constant $K(G)$  is independent of the selection of starting node $i$~\cite{LeLo02}.
It has found various applications in many fields. For example, it has been recently   applied to gauge the robotic surveillance efficiency~\cite{PaAgBu15}.

The hitting time  $T_{ij}$ for random walks on graph $G$ is encoded in the eigenvalues and eigenvectors of its normalized adjacency (or Laplacian) matrix.
\newtheorem{theorem}{Theorem}
\begin{theorem}\label{HitTime}\cite{Lo93}
For random walks on a simple connected graph $G$, the  hitting time $T_{ij}$ from one node $i$ to another node $j$ is
\begin{equation*}
  T_{ij}=2m\sum_{k=2}^n \frac{1}{1-\lambda_k}
  \left(\frac{v_{kj}^2}{d_j}-\frac{v_{ki}v_{kj}}{\sqrt{d_i d_j}}\right).  
\end{equation*}
In particular, when  $G$ is a bipartite graph with $V(G)=V_1 \cup V_2$, then
\begin{equation*}
T_{ij}= 2m \sum_{k=2}^{n-1} \frac{1}{1-\lambda_k}
     \left(\frac{v_{kj}^2}{d_j}-\frac{v_{ki}v_{kj}}{\sqrt{d_id_j}}\right),
\end{equation*}
if $i$ and $j$ are both in $V_1$ or $V_2$;
\begin{equation*}
T_{ij}= 2m \sum_{k=2}^{n-1} \frac{1}{1-\lambda_k}
      \left(\frac{v_{kj}^2}{d_j}-\frac{v_{ki}v_{kj}}{\sqrt{d_id_j}}\right)+1,
\end{equation*}
otherwise.
\end{theorem}

In contrast, the Kemeny's constant of $G$ is only dependent on the eigenvalues of $P$.
\begin{lemma}\label{lemmaKem}\cite{butler2016algebraic}
Let $G$ be a simple connected graph with $n$ nodes. Then
\begin{equation*}
K(G)=\sum_{j=1}^{n} \pi_j T_{ij}= \sum_{i=2}^n\frac{1}{1-\lambda_i},
\end{equation*}
where $1=\lambda_1>\lambda_2\geq\ldots\geq\lambda_n\geq-1$ are eigenvalues of matrix $P$.
\end{lemma}

\subsection{Electrical Networks}

For a simple connected graph $G$, we can define a corresponding electrical network $G^*$,  which is obtained from   $G$ by replacing each edge in  $G$ with  a unit resistor~\cite{DoSn84}. The resistance distance $r_{ij}$ between a pair of nodes $i$ and $j$ in $G$ is equal to the effective resistance between $i$ and $j$ in $G^*$.  Similar to the hitting time $T_{ij}$, resistance distance
$r_{ij}$ can also be expressed in terms of the eigenvalues and eigenvectors of normalized adjacency matrix $P$.
\begin{lemma}\cite{ChZh07}
Let $G$ be a simple connected  graph.  Then   resistance distance $r_{ij}$ between nodes $i$  and $j$ is
\begin{equation*}
    r_{ij}  =   \sum_{k=2}^n \frac{1}{1-\lambda_k}
                \left(\frac{v_{ki}}{\sqrt{d_i }}-\frac{v_{kj}}{\sqrt{d_j}}\right)^2.
\end{equation*}
\end{lemma}

\begin{lemma}\label{Foster}\cite{Fo49}
Let $G$ be a simple connected graph with $n$ nodes. Then the sum of resistance distances between all pairs of adjacent nodes in  $G$ is equivalent to $n-1$, i.e.
\begin{equation*}
  \sum_{ij\in E(G)}r_{ij}=n-1.
\end{equation*}
where the summation is taken over all the edges in $G$.
\end{lemma}

There are some intimate relationships between random walks on graphs and electrical networks. For example,
the resistance distance $r_{ij}$ is closely related to hitting times $T_{ij}$ and $T_{ji}$ of $G$, as stated in the  following lemma.
\begin{lemma}\cite{ChRaRuSm89}\label{lemmaRij}
For any pair of nodes $i$ and $j$ in a simple connected graph $G$ with $m$ edges, the following relation holds true:
\begin{equation*}
2mr_{ij}=T_{ij}+T_{ji}.
\end{equation*}
\end{lemma}

The resistance distance is an important quantity~\cite{GhBoSa08}. Various graph invariants based on resistance distances have been defined and studied. Among these invariants, the Kirchhoff index~\cite{KlRa93} is of vital importance.
\newtheorem{definition}{Definition}
\begin{definition}\label{defK}\cite{KlRa93}
The Kirchhoff index of a graph $G$ is defined as
\begin{equation*}
 \mathcal{K}(G)=\frac{1}{2}\sum_{i,j=1}^{n}r_{ij}=\sum_{\{i,j\}\subseteq V(G)}r_{ij}.
\end{equation*}
\end{definition}
Kirchhoff index  has found wide applications. For example, it can be used as  measures of the overall connectedness of a network~\cite{TiLe10}, the robustness of first-order consensus algorithm in noisy networks~\cite{PaBa14}, as well as the edge centrality of complex networks~\cite{LiZh18}.

In recent years,  several modifications for Kirchhoff index have been proposed, including additive degree-Kirchhoff index~\cite{GuFeYu12} and multiplicative degree-Kirchhoff index~\cite{ChZh07}. For a graph $G$, its additive degree-Kirchhoff index $\bar{\mathcal{K}}(G)$ and multiplicative degree-Kirchhoff index $\hat{\mathcal{K}}(G)$ are defined as
\begin{equation*}
  \bar{\mathcal{K}}(G)=\frac{1}{2}\sum_{i,j=1}^{n}(d_i+d_j)r_{ij}=\sum_{\{i,j\}\subseteq V(G)}(d_i+d_j)r_{ij}
\end{equation*}
and
\begin{equation*}
    \hat{\mathcal{K}}(G)=\frac{1}{2}\sum_{i,j=1}^{n}d_id_jr_{ij}=\sum_{\{i,j\}\subseteq V(G)}d_id_jr_{ij},
\end{equation*}
respectively.

It has been proved that $\hat{\mathcal{K}}(G)$ can be represented in terms of the eigenvalues of the matrix $P$.
\begin{lemma}\label{lemmaK*}\cite{ChZh07}
Let $G$ be a simple connected graph with $n$ nodes and $m$ edges. Then
\begin{equation*}
 \hat{\mathcal{K}}(G)=2m\sum_{i=2}^n\frac{1}{1-\lambda_i}.
\end{equation*}
\end{lemma}

\section{ $q$-triangulation Graphs and Their Matrices }

In this section, we introduce the  $q$-triangulation graph of a graph $G$,  which is a generalization of the traditional triangulation graph, since $1$-triangulation graph is exactly the triangulation graph. The triangulation of $G$, denoted by $R(G)$, is the graph obtained from $G$ by adding, for each edge $uv$ in $G$, a new node $x$ and two edges  $xu$ and $xv$. The triangulation graph can be easily extended to a general case.

\begin{definition}
Let $G$ be a simple connected graph. For a positive integer $q$, the $q$-triangulation graph of $G$, denoted by $R_q(G)$, is the graph obtained from $G$ by adding, for each edge $uv$ in $G$, $q$ disjoint paths of length $2$: $ux_1v$, $ux_2v$, $\ldots$, $ux_qv$.
\end{definition}

In what follows, for a quantity $Z$ of $G$, we use   $\tilde{Z}$  to denote the corresponding quantity associated with $R_q(G)$. Then it is easy to verify that in the $q$-triangulation graph $R_q(G)$, there are $\tilde{n}=n+mq$ nodes and $\tilde{m} = m(2q+1)$ edges.

Moreover,  the  node set $\tilde{V}:=V(R_q(G))$ of $ R_q(G)$ can be divided into two disjoint parts $V$ and $V^\prime$,  where $V$ is the set of old nodes inherited from $G$, while $V^\prime$ is the set of  new nodes generated in the process of performing $q$-triangulation operation on $G$. Moreover, $V^\prime$ can be further classified into $q$ parts as ${V}^\prime=V^{(1)} \cup  V^{(2)}\cup \cdots  \cup V^{(q)}$,  where each   $V^{(i)}$ ($i=1$, $2$,$\ldots$, $q$) contains $m$ new nodes produced by $m$ different edges in $G$. Namely,

\begin{equation} \label{div}
 \tilde{V}  =   V{\cup{V^{(1)}}} {\cup{V^{(2)}}} {\cup{...}} {\cup{V^{(q)}}}.
\end{equation}
By construction, for each old edge $uv$, there exists one and only one node $x$ in each $V^{(i)}$ ($i=1$, $2$,$\ldots$, $q$), satisfying  $\tilde{\Gamma} (x)=\{u,v\}$. Thus, for two different sets $V^{(i)}$ and $V^{(j)}$, the structural and dynamical properties of nodes belonging to them are equivalent to each other.

For $R_q(G)$,  its adjacency matrix $\tilde{A}$,  diagonal degree matrix $\tilde{D}$, and normalized adjacency matrix $\tilde{P}$, can be expressed in terms of related matrices of $G$ as

\begin{equation*}
\tilde{A}=
\left(
  \begin{array}{cccc}
    A           & B         & \cdots\   & B \\
    B^\top      & O         & \cdots\   & O  \\
    \vdots\     & \vdots\   & \ddots\   & \vdots\ \\
    B^\top      & O         &\cdots\    & O  \\
  \end{array}
\right),
\end{equation*}

\begin{equation*}
    \tilde{D}={\rm diag} \{(q+1)D,\underbrace{2I_m,...,2I_m}_q\},
\end{equation*}
and
\begin{equation}\label{MatrixP}
\begin{aligned}
&   \tilde{P}   =   \tilde{D}^{-\frac{1}{2}} \tilde{A} \tilde{D}^{-\frac{1}{2}}\\
=&  \frac{1}{\sqrt{2(q+1)}}
\left(
  \begin{array}{cccc}
    \sqrt{\frac{2}{q+1}}P    &D^{-\frac{1}{2}}B  & \cdots\   & D^{-\frac{1}{2}}B \\
    B^\top D^{-\frac{1}{2}} & O                 & \cdots\   & O \\
    \vdots\                 & \vdots\           & \ddots\   & O \\
    B^\top D^{-\frac{1}{2}} & O                 & \cdots\   & O \\
  \end{array}
\right),
\end{aligned}
\end{equation}
where $I_m$ is the $m \times m $ identity matrix.

\section{Eigenvalues and Eigenvectors of Normalized Adjacency Matrix
for\\ $q$-triangulation Graphs}

In this section, we study the eigenvalues and eigenvectors of normalized adjacency matrix $\tilde{P}$  for  $q$-triangulation graphs $R_q(G)$. We will show that both eigenvalues and eigenvectors for $\tilde{P}$ can be expressed in terms of those related quantities associated with graph $G$.

For the sake of convenience, for each eigenvalue $\lambda_i$ of $P$, we define $\Delta_i$ as $\Delta_i:= \lambda_i^2 + 2q(q+1)(1+\lambda_i)$. Thus, by Eqs.~(\ref{div}) and~(\ref{MatrixP}), we have the following result.

\begin{theorem}\label{Spectra}
Let $G$ be a simple connected graph with $n$ nodes and $m$ edges. Let $1=\lambda_1>\lambda_2\geq...\geq\lambda_n\geq-1$ be the eigenvalues of $P$, and let $v_1,v_2,...,v_n$ be their corresponding orthonormal eigenvectors. Then
\begin{enumerate}
  \item if $G$ is non-bipartite, then $\frac{\lambda_i\pm\sqrt{\Delta_i}}{2(q+1)}$ , $i=1$, $2$,$\ldots$, $n$ are eigenvalues of $\tilde{P}$, and the corresponding orthonormal eigenvectors are
      \begin{equation*}
        \sqrt{\frac{1}{2} \pm \frac{\lambda_i}{2\sqrt{\Delta_i}}}
        \left(
          \begin{array}{c}
            v_i \\
            \frac{\sqrt{2(q+1)}}{\lambda_i \pm \sqrt{\Delta_i}} B^\top D^{-\frac{1}{2}}v_i \\
            \vdots \\
            \frac{\sqrt{2(q+1)}}{\lambda_i \pm \sqrt{\Delta_i}} B^\top D^{-\frac{1}{2}}v_i \\
          \end{array}
        \right);
      \end{equation*}
      and $0$'s are eigenvalues of $\tilde{P}$ with multiplicity $mq-n$, with  their corresponding orthonormal eigenvectors being
      \begin{equation*}
             \\ \left(
                  \begin{array}{c}
                    0 \\
                    Y_z \\
                  \end{array}
                \right) ,~z=1,~2,\ldots,~mq-n,
      \end{equation*}
      where $Y_1,$ $Y_2$,$\ldots$, $Y_{mq-n}$ are an orthonormal basis of the kernel space of matrix
      \begin{equation*}
      C:={\underbrace{\left(
        \begin{array}{cccc}
            B & B & \cdots\ & B  \\
        \end{array}
      \right)}_q}.
      \end{equation*}
  \item if $G$ is bipartite, then $\frac{\lambda_i\pm\sqrt{\Delta_i}}{2(q+1)}$, $i=1, 2,\ldots, n-1$ are eigenvalues of $\tilde{P}$, and the corresponding orthonormal eigenvectors are
            \begin{equation*}
        \sqrt{\frac{1}{2} \pm \frac{\lambda_i}{2\sqrt{\Delta_i}}}
        \left(
          \begin{array}{c}
            v_i \\
            \frac{\sqrt{2(q+1)}}{\lambda_i \pm \sqrt{\Delta_i}} B^\top D^{-\frac{1}{2}}v_i \\
            \vdots \\
            \frac{\sqrt{2(q+1)}}{\lambda_i \pm \sqrt{\Delta_i}} B^\top D^{-\frac{1}{2}}v_i \\
          \end{array}
        \right);
      \end{equation*}
      $0$'s are eigenvalues of $\tilde{P}$ with multiplicity $mq-n+1$, with  their corresponding orthonormal eigenvectors being
      \begin{equation*}
             \\ \left(
                  \begin{array}{c}
                    0 \\
                    Y_z \\
                  \end{array}
                \right) , z=1, 2,\ldots, mq-n+1,
      \end{equation*}
      where $Y_1,Y_2,\ldots,Y_{mq-n+1}$ is an orthonormal basis of the kernel space of matrix
      \begin{equation*}
      C:={\underbrace{\left(
        \begin{array}{cccc}
            B & B & \cdots\ & B  \\
        \end{array}
      \right)}_q};
      \end{equation*}
      and $-\frac{1}{q+1}$ is an eigenvalue of $\tilde{P}$ of single degeneracy, with its corresponding eigenvector being
      \begin{equation*}
             \\ \left(
                  \begin{array}{c}
                    v_n \\
                    0 \\
                  \end{array}
                \right).
      \end{equation*}
\end{enumerate}
\end{theorem}

\begin{IEEEproof}
We first prove 1). Since $G$ is non-bipartite, by Lemma~\ref{MinEig}, every  eigenvalue $\lambda_i$ of $P$ is not equal to $-1$.  Notice that $Pv_i=\lambda_i v_i$. Then by Lemma~\ref{BBAD} and Eq.~(\ref{MatrixP}), it is easy to verify
\begin{equation*}
\begin{aligned}
        &\tilde{P}
        \left(
          \begin{array}{c}
            v_i \\
            \frac{\sqrt{2(q+1)}}{\lambda_i \pm \sqrt{\Delta_i}} B^\top D^{-\frac{1}{2}}v_i \\
            \vdots \\
            \frac{\sqrt{2(q+1)}}{\lambda_i \pm \sqrt{\Delta_i}} B^\top D^{-\frac{1}{2}}v_i \\
          \end{array}
        \right)
        =
        \left(
          \begin{array}{c}
            \frac{\lambda_i v_i}{q+1}+ \frac{q(1+\lambda_i)v_i}{\lambda_i\pm\sqrt{\Delta_i}} \\
            \frac{1}{\sqrt{2(q+1)}} B^\top D^{-\frac{1}{2}}v_i \\
            \vdots \\
            \frac{1}{\sqrt{2(q+1)}} B^\top D^{-\frac{1}{2}}v_i \\
          \end{array}
        \right)\\
        =&
        \frac{\lambda_i\pm\sqrt{\Delta_i}}{2(q+1)}
        \left(
          \begin{array}{c}
            v_i \\
            \frac{\sqrt{2(q+1)}}{\lambda_i \pm \sqrt{\Delta_i}} B^\top D^{-\frac{1}{2}}v_i \\
            \vdots \\
            \frac{\sqrt{2(q+1)}}{\lambda_i \pm \sqrt{\Delta_i}} B^\top D^{-\frac{1}{2}}v_i \\
          \end{array}
        \right),
\end{aligned}
\end{equation*}
which leads to our result through normalization.

For the zero eigenvalues, from Lemma~\ref{bRank}, ${\rm rank}(B)=n$ since $G$ is non-bipartite. Thus, ${\rm rank}(C)=n$ and $\dim({\rm Ker}(C))=mq-n$. Let $Y_1$, $Y_2$,$\ldots$, $Y_{mq-n}$ be an orthonormal basis of the kernel space of matrix $C$. It is easy to confirm  that
     $\left(
         \begin{array}{c}
                    0 \\
                    Y_z \\
         \end{array}
     \right)$, $z=1,2,\ldots, mq-n$, are eigenvectors for eigenvalues $0$ of matrix $\tilde{P}$. Moreover, together with the aforementioned eigenvectors, they constitute an orthonormal basis of $\tilde{P}$.

For 2), our proof is similar. We just need to verify that
\begin{equation*}
             \\ \tilde{P}\left(
                  \begin{array}{c}
                    v_n \\
                    0 \\
                  \end{array}
                \right)= -\frac{1}{q+1}
                \left(
                  \begin{array}{c}
                    v_n \\
                    0 \\
                  \end{array}
                \right),
\end{equation*}
which is trivial according to Eqs.~\eqref{eig=-1} and~\eqref{MatrixP}.
\end{IEEEproof}

Note that when $q=1$, Theorem~\ref{Spectra} reduces to the result in~\cite{XiZhCo16b}.

\section{Hitting Times for Random Walks on $q$-triangulation Graphs}

Theorem~\ref{Spectra} provides  complete information about the eigenvalues and eigenvectors of $\tilde{P}$ in terms of those of $P$. In this section, we use this information to determine two-node hitting time and Kemeny's constant for unbiased random walks on $R_q(G)$.

\subsection{Two-Node Hitting Time}
We first compute the hitting time from one node to another in $R_q(G)$. For this purpose, we  express the orthonormal eigenvectors of $R_q(G)$ in more explicit forms. By Eq.~\eqref{eig=1} and Theorem~\ref{Spectra}, we can easily derive the following results.
\begin{enumerate}
  \item The eigenvectors corresponding to eigenvalues $\frac{\lambda_1\pm\sqrt{\Delta_1}}{2(q+1)}=1, -\frac{q}{q+1}$ for matrix $\tilde{P}$ are
  \begin{small}
  \begin{equation}\label{re1A}
  \begin{aligned}
      \Big(&\sqrt{\frac{(q+1)d_1}{2m(2q+1)}}, \cdots, \sqrt{\frac{(q+1)d_n}{2m(2q+1)}},  \\&\frac{1}{\sqrt{m(2q+1)}},  \cdots, \frac{1}{\sqrt{m(2q+1)}} \Big)^\top
  \end{aligned}
  \end{equation}
  \end{small}
      and
  \begin{small}
  \begin{equation}\label{re1B}
  \begin{aligned}
          \Big(&\sqrt{\frac{qd_1}{2m(2q+1)}}, \cdots, \sqrt{\frac{qd_n}{2m(2q+1)}},  \\&-\sqrt{\frac{q+1}{mq(2q+1)}},  \cdots,    -\sqrt{\frac{q+1}{mq(2q+1)}} \Big)^\top
  \end{aligned}
  \end{equation}
  \end{small}
  respectively.
  \item If $G$ is non-bipartite, then $\frac{\lambda_i\pm\sqrt{\Delta_i}}{2(q+1)}$, $i=1$, $2$,$\ldots$, $n$ are eigenvalues of $\tilde{P}$, and the element of their orthonormal eigenvectors corresponding to node $j$ is
      \begin{small}
      \begin{equation}\label{re2}
         \left\{
            \begin{array}{ll}
                \sqrt{\frac{1}{2}\pm\frac{\lambda_i}{2\sqrt{\Delta_i}}}v_{ij}, &\hbox{$j\in V$;} \\
                \pm\sqrt{\frac{q+1}{\Delta_i\pm\lambda_i\sqrt{\Delta_i}}}\Big(\frac{v_{is}}{\sqrt{d_s}}+\frac{v_{it}}{\sqrt{d_t}}\Big), &\hbox{$j\in V^\prime$, $\tilde{\Gamma}(j)=\{s,t\}$.}
            \end{array}
        \right.
      \end{equation}
      \end{small}
  Moreover, for each $j\in V^\prime$ with $\tilde{\Gamma}(j)=\{s,t\}$,
  \begin{small}
  \begin{equation}\label{sum1}
    \sum_{z=1}^{mq-n} Y_{zj}^{2}=1-\frac{1}{mq}-\sum_{k=2}^n \frac{1}{(1+\lambda_k)q}\bigg(\frac{v_{ks}}{\sqrt{d_s}}+\frac{v_{kt}}{\sqrt{d_t}}\bigg)^2.
  \end{equation}
  \end{small}

  \item If $G$ is bipartite, then $\frac{\lambda_i\pm\sqrt{\Delta_i}}{2(q+1)}$, $i=1, 2,\ldots, n-1$ are eigenvalues of $\tilde{P}$, and the element of their orthonormal eigenvectors corresponding to node $j$ is
      \begin{small}
      \begin{equation}\label{re3}
         \left\{
            \begin{array}{ll}
                \sqrt{\frac{1}{2}\pm\frac{\lambda_i}{2\sqrt{\Delta_i}}}v_{ij}, &\hbox{$j\in V$;} \\
                \pm\sqrt{\frac{q+1}{\Delta_i\pm\lambda_i\sqrt{\Delta_i}}}\Big(\frac{v_{is}}{\sqrt{d_s}}+\frac{v_{it}}{\sqrt{d_t}}\Big), &\hbox{$j\in V^\prime$, $\tilde{\Gamma}(j)=\{s,t\}$.}
            \end{array}
        \right.
      \end{equation}
      \end{small}
  Moreover, for each $j\in V^\prime$ with $\tilde{\Gamma}(j)=\{s,t\}$,
  \begin{small}
  \begin{equation}\label{sum2}
    \sum_{z=1}^{mq-n+1} Y_{zj}^{2}=1-\frac{1}{mq}-\sum_{k=2}^{n-1} \frac{1}{(1+\lambda_k)q}\bigg(\frac{v_{ks}}{\sqrt{d_s}}+\frac{v_{kt}}{\sqrt{d_t}}\bigg)^2.
  \end{equation}
  \end{small}
\end{enumerate}

Now we present our results for hitting times of random walks on $R_q(G)$.
\begin{theorem}\label{qTriHT}
Let $G$ be a simple connected graph with $n$ nodes and $m$ edges. $R_q(G)$ is the $q$-triangulation graph of $G$ with $\tilde{V}=V\cup{V^\prime}$. Then
\begin{enumerate}
  \item if $i$, $j\in V$, then $\tilde{T}_{ij}=\frac{4q+2}{q+2} T_{ij}$;
  \item if $i\in V^\prime$, $j\in V$, $\tilde{\Gamma}(i)=\{s,t\}$, then
  \begin{equation*}
    \begin{aligned}
       \tilde{T}_{ij}=& 1+\frac{2q+1}{q+2}(T_{sj}+T_{tj});\\
       \tilde{T}_{ji}=& m(2q+1)-1\\
       &+\frac{2q+1}{2(q+2)}\big[2(T_{js}+T_{jt})-(T_{ts}+T_{st})\big];
    \end{aligned}
  \end{equation*}
  \item if $i$, $j\in V^\prime$, $j\in V$, $\tilde{\Gamma}(i)=\{s,t\}$, $\tilde{\Gamma}(j)=\{u,v\}$, then
    \begin{equation*}
    \begin{aligned}
       \tilde{T}_{ji}= &m(2q+1)+\frac{2q+1}{2(q+2)}\big[T_{su}+T_{tu}\\
       &+T_{sv}+T_{tv}-(T_{uv}+T_{vu})\big].
    \end{aligned}
    \end{equation*}
\end{enumerate}
\end{theorem}

\begin{IEEEproof}
Note that $\tilde{m}=m(2q+1)$, $\tilde{d}_i=(q+1)d_i$ if $i\in V$, and $\tilde{d}_i=2$ if $i\in V^\prime$.

We first prove 1). We distinguish two cases:  (a)   $G$ is non-bipartite, and   (b) $G$ is   bipartite.  
When $G$  is  non-bipartite,  by Theorems~\ref{HitTime} and~\ref{Spectra}, we have
\begin{small}
\begin{equation*}
    \begin{aligned}
         \tilde{T}_{ij} =&  2\tilde{m}\sum_{k=2}^{n}
                            \Bigg(\frac{1}{1-\frac{\lambda_k+\sqrt{\Delta_k}}{2(q+1)}}
                            \Big(\frac{1}{2}+\frac{\lambda_k}{2\sqrt{\Delta_k}}\Big)
                            +\frac{1}{1-\frac{\lambda_k-\sqrt{\Delta_k}}{2(q+1)}}\\
                         &  \Big(\frac{1}{2}-\frac{\lambda_k}{2\sqrt{\Delta_k}}\Big)\Bigg)
                            \bigg(\frac{v_{kj}^2}{(q+1)d_j}
                            -\frac{v_{ki}v_{kj}}{(q+1)\sqrt{d_id_j}}\bigg)\\
                        =&  2\tilde{m}\sum_{k=2}^{n}\frac{2q+2}{q+2}\frac{1}{1-\lambda_k}
                            \bigg(\frac{v_{kj}^2}{(q+1)d_j}-\frac{v_{kj}v_{ki}}{(q+1)\sqrt{d_id_j}}\bigg)\\
                       =&   \frac{4q+2}{q+2}\cdot2m\sum_{k=2}^{n}\frac{1}{1-\lambda_k}
                            \Big(\frac{v_{kj}^2}{d_j}-\frac{v_{kj}v_{ki}}{\sqrt{d_id_j}}\Big)=\frac{4q+2}{q+2}T_{ij}.
    \end{aligned}
\end{equation*}
\end{small}
When $G$ is  non-bipartite, the proof is similar.  Thus 1) is proved.


We continue to prove 2). Since $\tilde{\Gamma}(i)=\{s,t\}$,
\begin{equation*}
    \tilde{T}_{ij}=1+\frac{1}{2}\big(\tilde{T}_{sj}+\tilde{T}_{tj}\big)=1+\frac{2q+1}{q+2}(T_{sj}+T_{tj}).
\end{equation*}
While for $\tilde{T}_{ji}$, we also divide it into two cases:  (a)   $G$  is a non-bipartite graph, and   (b) $G$ is a  bipartite graph.  For the first case that  $G$ is non-bipartite,  by Theorems~\ref{HitTime} and~\ref{Spectra} and Eqs.~(\ref{re1A})-(\ref{sum1}), we have
\begin{small}
\begin{equation*}
  \begin{aligned}
    \tilde{T}_{ji}
    =&  2\tilde{m}\Bigg(\frac{1}{1+\frac{q}{q+1}}\frac{2q+1}{2mq(2q+1)}
        +\sum_{k=2}^n\frac{1}{2}
        \Big(\frac{v_{ks}}{\sqrt{d_s}}+\frac{v_{kt}}{\sqrt{d_t}}\Big)^2\\
    &   \bigg(\frac{1}{1-\frac{\lambda_k+\sqrt{\Delta_k}}{2(q+1)}}
        \frac{q+1}{\Delta_k+\lambda_k\sqrt{\Delta_k}}\\
    &   +\frac{1}{1-\frac{\lambda_k-\sqrt{\Delta_k}}{2(q+1)}}
        \frac{q+1}{\Delta_k-\lambda_k\sqrt{\Delta_k}}\bigg)\\
    &   -\sum_{k=2}^n\frac{v_{kj}}{\sqrt{2(q+1)d_j}}
        \Big(\frac{v_{ks}}{\sqrt{d_s}}+\frac{v_{kt}}{\sqrt{d_t}}\Big)\\
    &   \bigg(\frac{1}{1-\frac{\lambda_k+\sqrt{\Delta_k}}{2(q+1)}}
        \sqrt{\frac{q+1}{\Delta_k+\lambda_k\sqrt{\Delta_k}}}\\
    &   +\frac{1}{1-\frac{\lambda_k-\sqrt{\Delta_k}}{2(q+1)}}
        \sqrt{\frac{q+1}{\Delta_k-\lambda_k\sqrt{\Delta_k}}}\bigg)
        +\sum_{z=1}^{mq-n}\frac{Y_{zi}^2}{2} \Bigg)\\
    =&  m(2q+1)-1+\frac{2q+1}{q+2}
        2m\sum_{k=2}^n\frac{1}{1-\lambda_k}
        \bigg(\Big(\frac{v_{ks}^2}{d_s}-\frac{v_{ks}v_{kj}}{\sqrt{d_sd_j}}\Big)\\
    &   +\Big(\frac{v_{kt}^2}{d_t}-\frac{v_{kt}v_{kj}}{\sqrt{d_td_j}}\Big)
        -\frac{1}{2}\Big(\frac{v_{ks}}{\sqrt{d_s}}-\frac{v_{kt}}{\sqrt{d_t}}\Big)^2\bigg)\\
    =&  m(2q+1)-1+\frac{2q+1}{2(q+2)}\big[2(T_{js}+T_{jt})-(T_{ts}+T_{st})\big].
  \end{aligned}
\end{equation*}
\end{small}
If $G$ is bipartite, our proof is similar. 

We finally prove 3). Considering $\tilde{\Gamma}(i)=\{s,t\}$, $\tilde{\Gamma}(j)=\{u,v\}$,  we obtain
\begin{equation*}
\begin{aligned}
        \tilde{T}_{ij}
    =&   1+\frac{1}{2}(\tilde{T}_{sj}+\tilde{T}_{tj})\\
    =&   (2q+1)m+\frac{2q+1}{2(q+2)}\big[T_{su}+T_{tu}\\
     &   +T_{sv}+T_{tv}-(T_{uv}+T_{vu})\big].
\end{aligned}
\end{equation*}

This completes the proof.
\end{IEEEproof}

\subsection{Kemeny's Constant}

In addition to the two-node hitting time, the  Kemeny's constant of  $R_q(G)$ can also be expressed in terms of that of $G$.

\begin{theorem}\label{conKem}
Let $G$ be a simple connected graph with $n$ nodes and $m$ edges, and let $R_q(G)$ be the $q$-triangulation graph. Then
\begin{small}
\begin{equation*}
K(R_q(G))=\frac{4q+2}{q+2}K(G)+\frac{q^2+(4n-1)q+2n}{(q+2)(2q+1)}+mq-n.
\end{equation*}
\end{small}
\end{theorem}

\begin{IEEEproof}
Suppose that $1=\lambda_1>\lambda_2\geq...\geq\lambda_n\geq-1$ are eigenvalues of the matrix $P$. We first consider the case that $G$ is a non-bipartite graph. For this case, by Lemma~\ref{lemmaKem} and Theorem~\ref{Spectra}, we have
\begin{small}
\begin{equation*}
    \begin{aligned}
        K(R_q(G))   =   &   \sum_{k=2}^n \frac{1}{1-\frac{\lambda_k+\sqrt{\Delta_k}}{2(q+1)}}
                            +\sum_{k=2}^n
                            \frac{1}{1-\frac{\lambda_k-\sqrt{\Delta_k}}{2(q+1)}}\\
                        &   +\frac{1}{1+\frac{q}{q+1}}+mq-n\\
                    =   &   \sum_{k=2}^n
                            \big(\frac{2}{q+2}+\frac{4q+2}{q+2}\frac{1}{1-\lambda_k}\big)
                            +\frac{q+1}{2q+1}+mq-n\\
                    =   &   \frac{4q+2}{q+2}K(G)+\frac{q^2+(4n-1)q+2n}{(q+2)(2q+1)}+mq-n.
    \end{aligned}
\end{equation*}
\end{small}
For the other case that $G$ is bipartite, we can prove similarly.
\end{IEEEproof}

\section{Resistance Distances of $q$-triangulation Graphs }

In this section, we determine the two-node resistance distance, multiplicative degree-Kirchhoff index, additive degree-Kirchhoff index, and  Kirchhoff index of $R_q(G)$, in terms of those of $G$.

\subsection{Two-Node Resistance Distance }

We first determine  the resistance distance between any pair of nodes in $R_q(G)$.
\begin{theorem}\label{conReDis}
Let $G$ be a simple connected graph with $n$ nodes and $m$ edges, and let $R_q(G)$ be the $q$-triangulation graph with node set $\tilde{V}=V\cup{V^\prime}$. Then
\begin{enumerate}
  \item for $i,j \in V$,
    \begin{equation*}
      \tilde{r}_{ij}=\frac{2}{q+2}r_{ij};
    \end{equation*}
  \item for $i\in V^\prime$, $j\in V$ and $\tilde{\Gamma}(i)=\{s,t\}$,
    \begin{equation*}
    \tilde{r}_{ij}=\frac{1}{2}+\frac{2r_{sj}+2r_{tj}-r_{st}}{2(q+2)};
    \end{equation*}
  \item for $i,j \in V^\prime$, $\tilde{\Gamma}(i)=\{s,t\}$ and $\tilde{\Gamma}(j)=\{u,v\}$,
    \begin{equation*}
    \tilde{r}_{ij}=1+\frac{r_{su}+r_{tu}+r_{sv}+r_{tv}-r_{uv}-r_{st}}{2(q+2)}.
    \end{equation*}
\end{enumerate}
\end{theorem}
\begin{IEEEproof}
The results follow directly from Lemma~\ref{lemmaRij} and Theorem~\ref{qTriHT}.
\end{IEEEproof}

\subsection{Some Intermediary Results}

In the next subsections, we will derive the Kirchhoff index, the additive degree-Kirchhoff index and the multiplicative degree-Kirchhoff index for $R_q(G)$. In the computation of the first two graph invariants, we need  the following two properties for resistance distances in $R_q(G)$.

\begin{lemma}\label{afore1}
Let $G$ be a simple connected graph with $n$ nodes and $m$ edges, and let $R_q(G)$ be the $q$-triangulation graph of $G$ with node set $\tilde{V}=V\cup{V^\prime}$. Then
\begin{equation*}
  \sum_{i\in V^\prime}\sum_{j\in V}\tilde{r}_{ij}=
  \frac{q}{q+2}\bar{\mathcal{K}}(G)+\frac{mnq}{2}-\frac{n(n-1)q}{2(q+2)}.
\end{equation*}
\end{lemma}

\begin{IEEEproof}
Note that $\sum_{i\in V^\prime}\sum_{j\in V}\tilde{r}_{ij}$ can be divided into two sum terms as
\begin{equation}
    \sum_{i\in V^\prime}\sum_{j\in V}\tilde{r}_{ij}=
    \sum_{i\in V^\prime} \sum_{j\in \tilde{\Gamma}(i)} \tilde{r}_{ij}
    +\sum_{i\in V^\prime} \sum_{j\in V\backslash\tilde{\Gamma}(i)} \tilde{r}_{ij}.
\end{equation}
We next compute the above two sum terms separately.
\begin{enumerate}
  \item
  As for the first term, by Lemma~\ref{Foster}, we have
    \begin{equation}\label{Kf3}
        \begin{aligned}
            \sum_{i\in V^\prime} \sum_{j\in \tilde{\Gamma}(i)} \tilde{r}_{ij}
        =&  \sum_{ij\in \tilde{E}}\tilde{r}_{ij}
            -\sum_{ij\in E}\tilde{r}_{ij}\\
        =&  \Big(|\tilde{V}|-1\Big)
            -\frac{2}{q+2}\Big(|V|-1\Big)\\
        =&  mq+\frac{(n-1)q}{q+2}.
        \end{aligned}
    \end{equation}
  \item
  As for the second term, suppose that $\tilde{\Gamma}(i)=\{s,t\}$. According to Eq.~\eqref{div}, Lemma~\ref{Foster} and Theorem~\ref{conReDis}, we have
    \begin{equation}\label{Kf4}
        \begin{aligned}
            &   \sum_{i\in V^\prime} \sum_{j\in V\backslash\tilde{\Gamma}(i)} \tilde{r}_{ij}\\
            =&  \sum_{m=1}^{q}\sum_{i\in V^{(f)}} \sum_{j\in V\backslash\tilde{\Gamma}(i)}
                \Big(\frac{1}{2}+\frac{2r_{sj}+2r_{tj}-r_{st}}{2(q+2)}\Big)\\
            =&  q\sum_{i\in V^{(1)}} \sum_{j\in V\backslash\tilde{\Gamma}(i)}
                \Big(\frac{1}{2}+\frac{2r_{sj}+2r_{tj}-r_{st}}{2(q+2)}\Big)\\
            =&  \sum_{i\in V^{(1)}}\bigg(
                \frac{(n-2)q}{2}
                +\sum_{j\in V\backslash\tilde{\Gamma}(i)}\Big(\frac{q(r_{sj}+r_{tj})}{q+2}\\
             &  -\frac{(n-2)q}{2(q+2)}r_{st}\Big)\bigg).
        \end{aligned}
    \end{equation}

    For convenience, let $r_s$ be the sum of resistance distances between $s$ and all other nodes in graph $G$, that is,
    \begin{equation*}
        r_s=\sum_{j\in V \atop j\neq s} r_{sj}.
    \end{equation*}
    Thus, Eq.~\eqref{Kf4} can be rewritten as
    \begin{equation}\label{Kf5}
        \begin{aligned}
            &   \sum_{i\in V^\prime}
                \sum_{j\in V\backslash\tilde{\Gamma}(i)} \tilde{r}_{ij}\\
            =&  \sum_{i\in V^{(1)}}\bigg(
                \frac{(n-2)q}{2}
                +\frac{q(r_s+r_t)}{q+2}
                -\frac{(n+2)q}{2(q+2)}r_{st}\bigg)\\
            =&  \frac{m(n-2)q}{2}
                +\frac{q}{q+2}\sum_{i\in V^{(1)}}\big(r_s+r_t\big)\\
             &  -\frac{(n+2)q}{2(q+2)}\sum_{i\in V^{(1)}}r_{st}.
        \end{aligned}
    \end{equation}

    The term $\frac{q}{q+2}\sum_{i\in V^{(1)}}\big(r_s+r_t\big)$ can be further computed as
    \begin{equation}\label{Kf6}
    \begin{aligned}
         &  \frac{q}{q+2}\sum_{i\in V^{(1)}}\big(r_s+r_t\big)\\
        =&  \frac{q}{q+2}\sum_{st\in E}\big(r_s+r_t\big)
        =   \frac{q}{q+2}\sum_{s\in V}d_s r_s\\
        =&  \frac{q}{q+2}\sum_{\{i,j\}\subseteq V}(d_i+d_j)r_{ij}
        =   \frac{q}{q+2}\Bar{\mathcal{K}}(G).
    \end{aligned}
    \end{equation}
    Further, by Lemma~\ref{Foster}, the term $\frac{(n+2)q}{2(q+2)}\sum_{i\in V^{(1)}}r_{st}$ can be evaluated as
    \begin{equation}\label{Kf7}
    \begin{aligned}
        \frac{(n+2)q}{2(q+2)}\sum_{i\in V^{(1)}}r_{st}=&
        \frac{(n+2)q}{2(q+2)}\sum_{st\in E}r_{st}\\
        =&\frac{(n+2)(n-1)q}{2(q+2)}.
    \end{aligned}
    \end{equation}

    Plugging Eqs.~\eqref{Kf6} and~\eqref{Kf7} into Eq.~\eqref{Kf5} gives
    \begin{small}
    \begin{equation}\label{Kf8}
    \begin{aligned}
         &\sum_{i\in V^\prime} \sum_{j\in V\backslash\tilde{\Gamma}(i)} \tilde{r}_{ij}\\
         =& \frac{m(n-2)q}{2}+\frac{q}{q+2}\bar{\mathcal{K}}(G)
         -\frac{(n+2)(n-1)q}{2(q+2)}.
    \end{aligned}
    \end{equation}
    \end{small}
\end{enumerate}
Combining Eqs.~\eqref{Kf3} and~\eqref{Kf8} gives the desired result.
\end{IEEEproof}

\begin{lemma}\label{afore2}
Let $G$ be a connected graph with $n$ nodes and $m$ edges, and let $R_q(G)$ be the $q$-triangulation graph of $G$ with node set $\tilde{V}=V\cup{V^\prime}$. Then
\begin{small}
\begin{equation*}
  \sum_{\{i,j\}\subseteq V^\prime} \tilde{r}_{ij}
    =\frac{q^2}{2(q+2)}\hat{\mathcal{K}}(G)+\frac{mq(mq-1)}{2}-\frac{m(n-1)q^2}{2(q+2)}.
\end{equation*}
\end{small}
\end{lemma}

\begin{IEEEproof}
  Suppose that $\tilde{\Gamma}(i)=\{s,t\}$ and $\tilde{\Gamma}(j)=\{u,v\}$. Then by Theorem~\ref{conReDis}, we obtain
  \begin{equation}\label{A2_1}
    \begin{aligned}
        &   \sum_{\{i,j\}\subseteq V^\prime} \tilde{r}_{ij}\\
        =&    \sum_{\{i,j\}\subseteq V^\prime}
            \bigg(1+\frac{r_{su}+r_{tu}+r_{sv}+r_{tv}}{2(q+2)}-\frac{r_{st}+r_{uv}}{2(q+2)}\bigg)\\
        =&  \frac{mq(mq-1)}{2}
            +\sum_{\{i,j\}\subseteq V^\prime}\frac{r_{su}+r_{tu}+r_{sv}+r_{tv}}{2(q+2)}\\
         &  -\sum_{\{i,j\}\subseteq V^\prime}\frac{r_{st}+r_{uv}}{2(q+2)}.
    \end{aligned}
  \end{equation}

  We now compute the second term in Eq.~\eqref{A2_1}. It is not easy to evaluate it directly. Thus, we will compute it in an alternative way. For any pair of nodes $\{k,l\}\subseteq V$, we consider how many times $r_{kl}$ appears in the summation. Observe that $r_{kl}$ is summed once if and only if there exists a unique subset $\{i,j\}\subseteq V^\prime$ such that $k\in\tilde{\Gamma}(i)$ and $l\in\tilde{\Gamma}(j)$. Thus, our problem can be simplified and converted to the following one: how many pairwise different aforementioned subsets exist?
  It is not difficult to see that if $k$ is not adjacent to $l$, there exist $q^2d_kd_l$ subsets; and that if $kl\in E$, there exist $q^2d_kd_l-q$ such subsets. Thus, once again by Lemma~\ref{Foster}, we have
  \begin{equation}\label{A2_2}
    \begin{aligned}
         &  \sum_{\{i,j\}\subseteq V^\prime}\frac{r_{su}+r_{tu}+r_{sv}+r_{tv}}{2(q+2)}\\
        =&  \sum_{\{k,l\}\subseteq V \atop kl\notin E} \frac{q^2d_kd_l}{2(q+2)}r_{kl}
            +\sum_{\{k,l\}\subseteq V \atop kl\in E} \frac{q^2d_kd_l-q}{2(q+2)}r_{kl}\\
        =&  \sum_{\{k,l\}\subseteq V} \frac{q^2d_kd_l}{2(q+2)}r_{kl}
            -\frac{q}{2(q+2)}\sum_{\{k,l\}\subseteq V \atop kl\in E} r_{kl}\\
        =&  \frac{q^2}{2(q+2)}\sum_{\{k,l\}\subseteq V} d_kd_lr_{kl}-\frac{(n-1)q}{2(q+2)}\\
        =&  \frac{q^2}{2(q+2)}\hat{\mathcal{K}}(G)-\frac{(n-1)q}{2(q+2)}.
    \end{aligned}
  \end{equation}

  We proceed to evaluate the third term in Eq.~\eqref{A2_1}. Note that for any two different nodes $i$ and $j$ in $V^\prime$, if their neighbors are the same, i.e., $\tilde{\Gamma}(i)=\tilde{\Gamma}(j)=\{s,t\}$, we use $i\sim j$ to denote this relation. Otherwise, the sets of their neighbors are different, we call $i\nsim j$. According to these two relations and Eq.~\eqref{div}, it follows that
  \begin{small}
  \begin{equation}\label{A2_3}
    \begin{aligned}
        &   \sum_{\{i,j\}\subseteq V^\prime}\frac{r_{st}+r_{uv}}{2(q+2)}
        =   \frac{1}{4(q+2)}\sum_{i\in V^\prime}\sum_{j\in V^\prime}\big(r_{st}+r_{uv}\big)\\
        =&  \frac{1}{4(q+2)}\sum_{f=1}^{q}\sum_{i\in V^{(f)}}\Bigg(
            \sum_{i\nsim j} \big(r_{st}+r_{uv}\big)
            +\sum_{i\sim j \atop i\neq j}\big(r_{st}+r_{st}\big)\Bigg)\\
        =&  \frac{1}{4(q+2)}q\sum_{st\in E}\Bigg(
            q\sum_{uv\in E \atop uv\neq st} \big(r_{st}+r_{uv}\big)
            +2(q-1)r_{st}\Bigg)\\
        =&  \frac{q}{4(q+2)}\sum_{st\in E}\Bigg(
            q\sum_{uv\in E}r_{uv}
            +(mq-2)r_{st}\Bigg).
    \end{aligned}
  \end{equation}
  \end{small}
  By Lemma~\ref{Foster}, Eq.~\eqref{A2_3} can be recast as
  \begin{equation}\label{A2_4}
    \begin{aligned}
        &   \sum_{\{i,j\}\subseteq V^\prime}\frac{r_{st}+r_{uv}}{2(q+2)}\\
        =&  \frac{q}{4(q+2)}\sum_{st\in E}\Big((n-1)q+(mq-2)r_{st}\Big)\\
        =&  \frac{m(n-1)q^2}{4(q+2)}+\frac{(mq-2)(n-1)q}{4(q+2)}.
    \end{aligned}
  \end{equation}
  Plugging Eqs.~\eqref{A2_2} and~\eqref{A2_4} into Eq.~\eqref{A2_1} gives the result.
\end{IEEEproof}

\subsection{The multiplicative degree-Kirchhoff index}

We first  determine the multiplicative degree-Kirchhoff index for $R_q(G)$.
\begin{theorem}\label{conK*}
Let $G$ be a simple connected graph with $n$ nodes and $m$ edges, and let $R_q(G)$ be its $q$-triangulation graph. Then
\begin{small}
\begin{equation*}
\begin{aligned}
    \hat{\mathcal{K}}(R_q(G))
=&  \frac{2(2q+1)^2}{q+2}\hat{\mathcal{K}}(G)
    +2m\bigg(\frac{q^2+(4n-1)q+2n}{q+2}\\
&   +(mq-n)(2q+1)\bigg).
\end{aligned}
\end{equation*}
\end{small}
\end{theorem}

\begin{IEEEproof}
According to Lemmas~\ref{lemmaKem} and \ref{lemmaK*}, Theorem~\ref{conK*} is an obvious consequence of Theorem~\ref{conKem}.
\end{IEEEproof}

\subsection{The addictive degree-Kirchhoff index}

We continue to determine the additive degree-Kirchhoff index for $R_q(G)$.
\begin{theorem}\label{con6}
Let $G$ be a connected graph with $n$ nodes and $m$ edges, and let $R_q(G)$ be the $q$-triangulation graph. Then
\begin{equation*}
\begin{aligned}
    \bar{\mathcal{K}}(R_q(G))
=&  \frac{2(2q+1)}{q+2}\bar{\mathcal{K}}(G)
    +\frac{2q(2q+1)}{q+2}\hat{\mathcal{K}}(G)\\
&   +m^2q(3q+1) -  mq(2n-1)\\
&   +\frac{(5m-n)(n-1)q}{q+2}.
\end{aligned}
\end{equation*}
\end{theorem}

\begin{IEEEproof}
By definition of the addictive degree-Kirchhoff index, we have
\begin{equation}\label{Kf+1}
    \begin{aligned}
    \bar{\mathcal{K}}(R_q(G))=&  \sum_{\{i,j\}\subseteq V\cup V^\prime} (\tilde{d}_i+\tilde{d}_j) \tilde{r}_{ij}\\
                    =&  \sum_{\{i,j\}\subseteq V} (\tilde{d}_i+\tilde{d}_j) \tilde{r}_{ij}
                        +\sum_{i\in V^\prime}\sum_{j\in V}(\tilde{d}_i+\tilde{d}_j) \tilde{r}_{ij}\\
                     &  +\sum_{\{i,j\}\subseteq V^\prime} (\tilde{d}_i+\tilde{d}_j) \tilde{r}_{ij}.
    \end{aligned}
\end{equation}
We now compute the three sum terms on the last row of Eq.~\eqref{Kf+1} one by one.

For the first sum term, by Theorem~\ref{conReDis}, we have
\begin{small}
\begin{equation}\label{Kf+2}
\begin{aligned}
    &   \sum_{\{i,j\}\subseteq V} (\tilde{d}_i+\tilde{d}_j) \tilde{r}_{ij}\\
    =&  \sum_{\{i,j\}\subseteq V} (q+1)(d_i+d_j) \frac{2}{q+2}r_{ij}
    =   \frac{2(q+1)}{q+2}\bar{\mathcal{K}}(G).
\end{aligned}
\end{equation}
\end{small}

For the second sum term, it can be evaluated as
\begin{equation}\label{Kf+3}
\begin{aligned}
     &  \sum_{i\in V^\prime}\sum_{j\in V}(\tilde{d}_i+\tilde{d}_j) \tilde{r}_{ij}\\
    =&  \sum_{i\in V^\prime}\sum_{j\in V}(2+(q+1)d_j) \tilde{r}_{ij}\\
    =&  2\sum_{i\in V^\prime}\sum_{j\in V}\tilde{r}_{ij}
        +(q+1)\sum_{i\in V^\prime}\sum_{j\in V}d_j\tilde{r}_{ij}.
\end{aligned}
\end{equation}
By Lemma~\ref{afore1}, we have
\begin{equation}\label{Kf+4}
    2\sum_{i\in V^\prime}\sum_{j\in V}\tilde{r}_{ij}=
    \frac{2q}{q+2}\bar{\mathcal{K}}(G)+mnq-\frac{n(n-1)q}{q+2}.
\end{equation}
On the other hand, by Lemma~\ref{Foster} and Theorem~\ref{conReDis},
\begin{equation}\label{Kf+5}
\begin{aligned}
        &     (q+1)\sum_{i\in V^\prime}\sum_{j\in V}d_j\tilde{r}_{ij}\\
        =&    (q+1)\sum_{i\in V^\prime}\sum_{j\in V}d_j
              \bigg(\frac{1}{2}+\frac{2r_{sj}+2r_{tj}-r_{st}}{2(q+2)}\bigg)\\
        =&    \frac{q+1}{2}\sum_{i\in V^\prime}\sum_{j\in V}d_j
              +\frac{q+1}{q+2}\sum_{i\in V^\prime}\sum_{j\in V}d_j\Big(r_{sj}+r_{tj}\Big)\\
         &    -\frac{q+1}{2(q+2)}\sum_{i\in V^\prime}\sum_{j\in V}d_jr_{st}\\
        =&    \frac{q+1}{2}\sum_{i\in V^\prime}2m
              +\frac{q+1}{q+2}\sum_{i\in V^\prime}\sum_{j\in V}d_j\Big(r_{sj}+r_{tj}\Big)\\
         &    -\frac{q+1}{2(q+2)}\sum_{i\in V^\prime}2mr_{st}\\
        =&    m^2q(q+1)
              +\frac{q+1}{q+2}\sum_{i\in V^\prime}\sum_{j\in V}d_j
              \Big(r_{sj}+r_{tj}\Big)\\
         &    -\frac{m(n-1)q(q+1)}{q+2}.
\end{aligned}
\end{equation}
For the middle part of Eq.~\eqref{Kf+5}, we obtain
\begin{equation}\label{Kf+6}
  \begin{aligned}
      &   \frac{q+1}{q+2}\sum_{i\in V^\prime}\sum_{j\in V}d_j
          \big(r_{sj}+r_{tj}\big)\\
      =&  \frac{q+1}{q+2} q\sum_{i\in V^{(1)}}\sum_{j\in V}d_j
          \big(r_{sj}+r_{tj}\big)\\
      =&  \frac{q(q+1)}{q+2}\sum_{j\in V}\sum_{i\in V^{(1)}}d_j\big(r_{sj}+r_{tj}\big)\\
      =&  \frac{q(q+1)}{q+2}\sum_{j\in V}\sum_{k\in V}d_j d_k r_{kj}
      =   \frac{2q(q+1)}{q+2}\hat{\mathcal{K}}(G).
  \end{aligned}
\end{equation}
Combining Eqs.~\eqref{Kf+3}-\eqref{Kf+6} yields
\begin{equation}\label{Kf+7}
  \begin{aligned}
   &    \sum_{i\in V^\prime}\sum_{j\in V}
        (\tilde{d}_i+\tilde{d}_j) \tilde{r}_{ij}\\
  =&    \frac{2q}{q+2}\bar{\mathcal{K}}(G)
        +\frac{2q(q+1)}{q+2}\hat{\mathcal{K}}(G)+m^2q^2+m^2q\\
  &     +mnq-\frac{(mq+m+n)(n-1)q}{q+2}.
  \end{aligned}
\end{equation}

For the third sum term in Eq.~\eqref{Kf+1}, by Lemma~\ref{afore2}, we have
\begin{small}
\begin{equation}\label{Kf+8}
\begin{aligned}
&   \sum_{\{i,j\}\subseteq V^\prime} (\tilde{d}_i+\tilde{d}_j) \tilde{r}_{ij}\\
=&  4\bigg(\frac{q^2}{2(q+2)}
    \hat{\mathcal{K}}(G)+\frac{mq(mq-1)}{2}
    -\frac{m(n-1)q^2}{2(q+2)}\bigg)\\
=&  \frac{2q^2}{q+2}\hat{\mathcal{K}}(G)+2mq(mq-1)
    -\frac{2m(n-1)q^2}{q+2}.
\end{aligned}
\end{equation}
\end{small}

Substituting Eqs.~\eqref{Kf+2}, \eqref{Kf+7} and~\eqref{Kf+8} back into Eq.~\eqref{Kf+1}, our proof is completed after simple calculations.
\end{IEEEproof}

\subsection{The Kirchhoff index}

We finally determine the Kirchhoff index for $R_q(G)$.
\begin{theorem}\label{conKir}
Let $G$ be a connected graph with $n$ nodes, and let $R_q(G)$ be the $q$-triangulation graph. Then
\begin{small}
\begin{equation*}
\begin{aligned}
    \mathcal{K}(R_q(G))  =&  \frac{2}{q+2}\mathcal{K}(G)+\frac{q}{q+2}\bar{\mathcal{K}}(G)+\frac{q^2}{2(q+2)}\hat{\mathcal{K}}(G)\\
                    &+  \frac{m^2q^2}{2}+\frac{(2m-n)(n-1)q}{2(q+2)}.
\end{aligned}
\end{equation*}
\end{small}
\end{theorem}

\begin{IEEEproof}
According to Definition~\ref{defK} and Eq.~\eqref{div}, we have
\begin{equation}\label{Kf1}
  \begin{aligned}
    \mathcal{K}(R_q(G))  =&  \sum_{\{i,j\}\subseteq\tilde{V}} \tilde{r}_{ij}
                         =   \sum_{\{i,j\}\subseteq V\cup V^\prime} \tilde{r}_{ij}\\
                         =&  \sum_{\{i,j\}\subseteq V} \tilde{r}_{ij}
                             +\sum_{i\in V^\prime} \sum_{j\in V} \tilde{r}_{ij}
                             +\sum_{\{i,j\}\subseteq V^\prime} \tilde{r}_{ij}.
  \end{aligned}
\end{equation}
Below we shall compute the three sum terms in Eq.~\eqref{Kf1} separately.

For the first sum term, by Theorem~\ref{conReDis},
\begin{equation}\label{Kf2}
    \sum_{\{i,j\}\subseteq V} \tilde{r}_{ij}=\sum_{\{i,j\}\subseteq V} \frac{2}{q+2} r_{ij}
    =\frac{2}{q+2}\mathcal{K}(G).
\end{equation}

For the second sum term, by Lemma~\ref{afore1}, we obtain
\begin{equation}\label{Kf9}
  \sum_{i\in V^\prime}\sum_{j\in V}\tilde{r}_{ij}=
  \frac{q}{q+2}\bar{\mathcal{K}}(G)+\frac{mnq}{2}-\frac{n(n-1)q}{2(q+2)}.
\end{equation}

For the third sum term, by Lemma~\ref{afore2}, we have
\begin{small}
\begin{equation}\label{Kf10}
  \sum_{\{i,j\}\subseteq V^\prime} \tilde{r}_{ij}
  =\frac{q^2}{2(q+2)}\hat{\mathcal{K}}(G)+\frac{mq(mq-1)}{2}-\frac{m(n-1)q^2}{2(q+2)}.
\end{equation}
\end{small}

Plugging Eqs.~\eqref{Kf2}-\eqref{Kf10} back into Eq.~\eqref{Kf1} leads to the desired result.
\end{IEEEproof}

\section{Properties of Iterated $q$-triangulation graphs and Their Applications}

The $q$-triangulation graphs have found many applications in physics and network science. For example, by iteratively applying   $q$-triangulation operation to $3$-clique, a complete graph with $3$ nodes and $3$ edges, we can obtain a family of scale-free small-world networks, called  pseudofractal scale-free webs~\cite{DoGoMe02,ZhRoZh07}, which have attracted considerable attention~\cite{Ya06,XiZhCo16b,ShLiZh17,ShLiZh18}. In this section, we study the properties of iterated $q$-triangulation graphs, based on which we further obtain exact expressions for some interesting quantities for  pseudofractal scale-free webs.

\subsection{Definition  of Iterated $q$-triangulation Graphs}

The family of iterated $q$-triangulation graphs $R_{q,k}(G)$ of a graph $G$ is defined as follows. For $k=0$, $R_{q,0}(G)=G$. For $k\geq 1$, $R_{q,k}(G)$ is obtained from $R_{q, k-1}(G)$ by performing the $q$-triangulation operation on $R_{q, k-1}(G)$. In other words, $R_{q,k}(G)=R_q(R_{q,k-1}(G))$.  
For a quantity $Z$ of $G$, we use  $Z_{q,k}$  to denote the corresponding quantity associated with $R_{q,k}(G)$.  Then, in $R_{q,k}(G)$, the number of edges is
\begin{equation}\label{iteratedm}
    m_{q,k} =   (2q+1) m_{q,k-1}=(2q+1)^{k}m,
\end{equation}
and the number of   nodes is
\begin{equation}\label{iteratedn}
    n_{q,k} =   n_{q,k-1}+q m_{q,k-1}
            =   \frac{m\big[(2q+1)^{k}-1\big]}{2}+n.
\end{equation}

\subsection{Formulas of Quantities for Iterated $q$-triangulation Graphs}

We here present expressions for some interesting quantities  for iterated $q$-triangulation graphs $R_{q,k}(G)$.

\subsubsection{Kemeny's constant}
\begin{theorem}\label{iteratedKem}
Let $G$ be a connected graph with $n$ nodes and $m$ edges. Then
\begin{small}
\begin{equation*}
    \begin{aligned}
    K_{q,k}   =&  \Big(\frac{4q+2}{q+2}\Big)^kK_{q,0}
                +\frac{m(2q+3)}{2(2q+1)}
                \Big[(2q+1)^k-\Big(\frac{4q+2}{q+2}\Big)^k\Big]\\
            &   +\bigg(\frac{(q-1)}{3(2q+1)}
                +\frac{m-2n}{6}\bigg)
                \Big[\Big(\frac{4q+2}{q+2}\Big)^k-1\Big].
    \end{aligned}
\end{equation*}
\end{small}
\end{theorem}

\begin{IEEEproof}
According to Theorem~\ref{conKem} and Eqs.~\eqref{iteratedm} and~\eqref{iteratedn}, we have
\begin{equation*}
    \begin{aligned}
    K_{q,k} =&  \frac{4q+2}{q+2}K_{q,k-1}
                +\frac{q^2+(4n_{q,k-1}-1)q+2n_{q,k-1}}{(q+2)(2q+1)}\\
             &  +m_{q,k-1}q-n_{q,k-1}\\
            =&  \frac{4q+2}{q+2}K_{q,k-1}
                +\frac{mq(2q+3)(2q+1)^{k-1}}{2(q+2)}\\
             &  +\frac{q(q-1)}{(q+2)(2q+1)}
                +\frac{(m-2n)q}{2(q+2)}.
    \end{aligned}
\end{equation*}
Dividing  both sides by  $\Big(\frac{4q+2}{q+2}\Big)^k$, we obtain
\begin{equation*}
\begin{aligned}
    &   \Big(\frac{q+2}{4q+2}\Big)^k K_{q,k}
        -\Big(\frac{q+2}{4q+2}\Big)^{k-1} K_{q,k-1}\\
    =&  \frac{mq(2q+3)}{2(q+2)(2q+1)}\Big(\frac{q+2}{2}\Big)^{k}\\
    &   +\bigg(\frac{q(q-1)}{(q+2)(2q+1)}
        +\frac{(m-2n)q}{2(q+2)}\bigg)\Big(\frac{q+2}{4q+2}\Big)^k.
\end{aligned}
\end{equation*}
By properties of geometric sequences, we have
\begin{equation*}
\begin{aligned}
    &   \Big(\frac{q+2}{4q+2}\Big)^k K_{q,k}
        -\Big(\frac{q+2}{4q+2}\Big)^{0} K_{q,0}\\
    =&  \frac{m(2q+3)}{2(2q+1)}\Big[\Big(\frac{q+2}{2}\Big)^{k}-1\Big]\\
     &  +\bigg(\frac{(q-1)}{3(2q+1)}
        +\frac{m-2n}{6}\bigg)\Big[1-\Big(\frac{q+2}{4q+2}\Big)^k\Big],
\end{aligned}
\end{equation*}
which leads to the result through simple calculations.
\end{IEEEproof}

\subsubsection{The multiplicative degree-Kirchhoff index}
\begin{theorem}\label{iteratedK*}
Let $G$ be a connected graph with $n$ nodes and $m$ edges. Then
\begin{small}
\begin{equation*}
\begin{aligned}
        \hat{\mathcal{K}}_{q,k}
    =&  \Big(\frac{2(2q+1)^2}{q+2}\Big)^k\Hat{\mathcal{K}}_{q,0}
        +\frac{m^2(2q+3)}{2q+1}
        \Big[(2q+1)^{2k}\\
    &   -\Big(\frac{2(2q+1)^2}{q+2}\Big)^k\Big]
        +\bigg(\frac{2m(q-1)}{3(2q+1)}
        +\frac{m(m-2n)}{3}\bigg)\\
    &   \Big[\Big(\frac{2(2q+1)^2}{q+2}\Big)^k-(2q+1)^k\Big].
\end{aligned}
\end{equation*}
\end{small}
\end{theorem}

\begin{IEEEproof}
By Lemmas~\ref{lemmaKem} and~\ref{lemmaK*}, the result follows directly from Theorem~\ref{iteratedKem}.
\end{IEEEproof}

\subsubsection{The addictive degree-Kirchhoff index}
\begin{theorem}\label{iteratedK+}
Let $G$ be a connected graph with $n$ nodes and $m$ edges. Then
\begin{small}
\begin{equation*}
\begin{aligned}
    \bar{\mathcal{K}}_{q,k}
    =&  \Big(\frac{2(2q+1)}{q+2}\Big)^k\bar{\mathcal{K}}_{q,0}
        +\Big[\Big(\frac{2(2q+1)^{2}}{q+2}\Big)^k-
        \Big(\frac{2(2q+1)}{q+2}\Big)^k\Big]\\
    &   \bigg(\frac{\hat{\mathcal{K}}_{q,0}}{2}
        -\frac{2(q+2)m^2+(2q+1)mn-m(q-1)}{3(2q+1)}\bigg)\\
    &   +\Big[(2q+1)^{2k}-\Big(\frac{2(2q+1)}{q+2}\Big)^k\Big]
        \frac{m^2(2q+3)(6q+11)}{4(2q+1)(2q+5)}\\
    &   +\Big[\Big(\frac{2(2q+1)}{q+2}\Big)^k-(2q+1)^k\Big]\\
    &   \Big(\frac{m}{2(2q+1)}+\frac{(q+2)m(m-2n+1)}{3(2q+1)}\Big)\\
    &   -\Big[\Big(\frac{2(2q+1)}{q+2}\Big)^k-1\Big]
        \frac{(m-2n)(m-2n+2)}{12}.
\end{aligned}
\end{equation*}
\end{small}
\end{theorem}

\begin{IEEEproof}
By Theorem~\ref{con6} and Eqs.~\eqref{iteratedm} and~\eqref{iteratedn}, we obtain
\begin{equation}\label{add}
\begin{aligned}
    \Bar{\mathcal{K}}_{q,k}
    =&  \frac{2(2q+1)}{q+2}\bar{\mathcal{K}}_{q,k-1}
        +\frac{2q(2q+1)}{q+2}\hat{\mathcal{K}}_{q,k-1}\\
    &   +(m_{q,k-1})^2q(3q+1)-m_{q,k-1}q(2n_{q,k-1}-1)\\
    &   +\frac{(5m_{q,k-1}-n_{q,k-1})(n_{q,k-1}-1)q}{q+2}\\
    =&  \frac{2(2q+1)}{q+2}\bar{\mathcal{K}}_{q,k-1}
        +\frac{2q(2q+1)}{q+2}\hat{\mathcal{K}}_{q,k-1}\\
    &   +(2q+1)^{2k-2}\frac{3m^2q(4q^2+8q+3)}{4(q+2)}\\
    &   +(2q+1)^{k-1}\frac{mq}{2(q+2)}\Big(2q(m-2n+1)-5\Big)\\
    &   -\frac{q}{4(q+2)}(m-2n)(m-2n+2).
\end{aligned}
\end{equation}
Inserting Theorem~\ref{iteratedK*} into Eq.~\eqref{add} gives
\begin{small}
\begin{equation}
\begin{aligned}
    \Bar{\mathcal{K}}_{q,k}    =&  \frac{2(2q+1)}{q+2}\bar{\mathcal{K}}_{q,k-1}
                        +\Big(\frac{2(2q+1)^2}{q+2}\Big)^{k-1}
                        \Big[\frac{2q(2q+1)}{q+2}\hat{\mathcal{K}}_{q,0}\\
                    &   -\frac{4mq}{3(q+2)}\Big(2(q+2)m+(2q+1)n-(q-1)\Big)\Big]\\
                    &   +(2q+1)^{2k-2}\frac{m^2q(2q+3)(6q+11)}{4(q+2)}\\
                    &   -(2q+1)^{k-1}
                        \bigg(\frac{mq}{2(q+2)}+\frac{mq(m-2n+1)}{3}\bigg)\\
                    &   -\frac{q}{4(q+2)}(m-2n)(m-2n+2).
\end{aligned}
\end{equation}
\end{small}
Dividing both sides by $\Big(\frac{2(2q+1)}{q+2}\Big)^k$, we obtain a geometric sequence, which is solved to yield the result.
\end{IEEEproof}
\subsubsection{The Kirchhoff index}
\begin{theorem}\label{iteratedKir}
Let $G$ be a connected graph with $n$ nodes and $m$ edges. Then
\begin{small}
\begin{equation*}
  \begin{aligned}
    &   \mathcal{K}_{q,k}\\
    =&  \Big(\frac{2}{q+2}\Big)^k\mathcal{K}_{q,0}
        +\Big[\Big(\frac{2(2q+1)^2}{q+2}\Big)^k
        -\Big(\frac{2}{q+2}\Big)^k\Big]\\
     &  \Big(\frac{\hat{\mathcal{K}}_{q,0}}{16}
        -\frac{m^2(q+2)}{12(2q+1)}-\frac{mn}{24}
        +\frac{m(q-1)}{24(2q+1)}\Big)\\
     &  +\Big[\Big(\frac{2(2q+1)}{q+2}\Big)^k
        -\Big(\frac{2}{q+2}\Big)^k\Big]
        \bigg(\frac{{\bar{\mathcal{K}}}_{q,0}}{4}
        -\frac{{\hat{\mathcal{K}}}_{q,0}}{8}\\
    &   -\frac{m^2(q+2)(2q-1)}{6(2q+1)(2q+5)}
        +\frac{m(2n(q-1)-q+4)}{12(2q+1)}
        -\frac{n(n-1)}{12}\bigg)\\
    &   +\Big[(2q+1)^{2k}-\Big(\frac{2}{q+2}\Big)^k\Big]
        \frac{m^2(2q+3)^2}{8(2q+1)(2q+5)}\\
    &   -\Big[(2q+1)^{k}-\Big(\frac{2}{q+2}\Big)^k\Big]\\
    &   \bigg(\frac{m(4q^2+12q+11)(m-2n)}{12(2q+1)(2q+5)}
        +\frac{m(4q^2+18q+23)}{12(2q+1)(2q+5)}\bigg)\\
    &   +\Big[\Big(\frac{2}{q+2}\Big)^k-1\Big]
        \frac{(m-2n)(m-2n+2)}{24}.\\
  \end{aligned}
\end{equation*}
\end{small}
\end{theorem}

\begin{IEEEproof}
By Theorems~\ref{conKir} and Eqs.~\eqref{iteratedm} and~\eqref{iteratedn}, we obtain
\begin{equation}\label{kir}
\begin{aligned}
    \mathcal{K}_{q,k}
    =&  \frac{2}{q+2}\mathcal{K}_{q,k-1}+\frac{q}{q+2}{\bar{\mathcal{K}}}_{q,k-1}
        +\frac{q^2}{2(q+2)}{\hat{\mathcal{K}}}_{q,k-1}  \\
    &   +\frac{(m_{q,k-1})^2q^2}{2}
        +\frac{(2m_{q,k-1}-n_{q,k-1})(n_{q,k-1}-1)q}{2(q+2)}\\
    =&  \frac{2}{q+2}\mathcal{K}_{q,k-1}+\frac{q}{q+2}{\bar{\mathcal{K}}}_{q,k-1}
        +\frac{q^2}{2(q+2)}{\hat{\mathcal{K}}}_{q,k-1}\\
    &   +(2q+1)^{2k-2}\frac{m^2q(4q^2+8q+3)}{8(q+2)}\\
    &   -(2q+1)^{k-1}\frac{mq(m-2n+3)}{4(q+2)}\\
    &   -\frac{q}{8(q+2)}(m-2n)(m-2n+2).
\end{aligned}
\end{equation}
According to Theorems~\ref{iteratedK*} and~\ref{iteratedK+}, Eq.~\eqref{kir} can be rewritten as
\begin{small}
\begin{equation}
\begin{aligned}
    \mathcal{K}_{q,k}
=&  \frac{2}{q+2}\mathcal{K}_{q,k-1}+
    \Big(\frac{2(2q+1)^2}{q+2}\Big)^{k-1}
    \bigg(\frac{q(q+1)}{2(q+2)}{\hat{\mathcal{K}}}_{q,0}\\
&   -\frac{2m^2q(q+1)}{3(2q+1)}
    -\frac{mnq(q+1)}{3(q+2)}
    +\frac{mq(q+1)(q-1)}{3(q+2)(2q+1)}\bigg)\\
&   +\Big(\frac{2(2q+1)}{q+2}\Big)^{k-1}
    \bigg(\frac{q}{q+2}{\bar{\mathcal{K}}}_{q,0}
    -\frac{q}{2(q+2)}{\hat{\mathcal{K}}}_{q,0}\\
&   -\frac{2m^2q(2q-1)}{3(2q+1)(2q+5)}
    +\frac{mq(2n(q-1)-q+4)}{3(q+2)(2q+1)}\\
&   -\frac{n(n-1)q}{3(q+2)}\bigg)
    +(2q+1)^{2k-2}
    \frac{m^2q(2q+3)^4}{8(q+2)(2q+1)(2q+5)}\\
&   -(2q+1)^{k-1}
    \bigg(\frac{mq(4q^2+12q+11)(m-2n)}{12(q+2)(2q+1)}\\
&   +\frac{mq(4q^2+18q+23)}{12(q+2)(2q+1)}\bigg)\\
&   -\frac{q}{24(q+2)}(m-2n)(m-2n+2).
\end{aligned}
\end{equation}
\end{small}
Dividing both sides by $\Big(\frac{2}{q+2}\Big)^k$ , we  derive our result through simple calculations.
\end{IEEEproof}

Our results in this section generalize those previously obtained for triangulation graphs~\cite{YaKl15}, but our computation method is much easier.

\subsection{Applications to pseudofractal scale-free webs}

\begin{figure}
\begin{center}
\includegraphics[width=7.0cm]{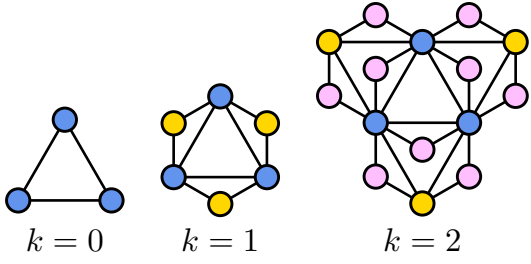}
\caption{Illustration of   pseudofractal scale-free webs $N_{1,0}$, $N_{1,1}$, and $N_{1,2}$.}\label{webq2}
\end{center}
\end{figure}

The family of pseudofractal scale-free webs~\cite{ZhRoZh07} is a particular example of iterated $q$-triangulation graphs. They are constructed in an iterative way. Let $N_{q,k}$, $q \geq 1$ and $k\geq 0$, denote the pseudofractal scale-free webs after $k$ iterations. For $k=0$, $N_{q,0}$ is the 3-node complete graph. For $k\geq1$, $N_{q,k}$ is obtained from $N_{q,k-1}$ by performing the $q$-triangulation operation on $N_{q,k-1}$. Thus, the pseudofractal scale-free webs are actually iterated $q$-triangulation graphs  $R_{q,k}(G)$ when $G$ is a $3$-clique. Fig.~\ref{webq2} illustrates the first several iterations for pseudofractal scale-free webs for  a particular case of $q=1$.

In the sequel, we provide some properties of the pseudofractal scale-free webs,   using the results obtained in last subsections.

For $N_{q,0}$, its adjacency matrix, diagonal degree matrix, and normalized adjacency matrix are
\begin{equation*}
    A(N_{q,0})=
         \left(
           \begin{array}{ccc}
             0      & 1         & 1 \\
             1      & 0         & 1 \\
             1      & 1         & 0 \\
           \end{array}
        \right),
\end{equation*}
\begin{equation*}
    D(N_{q,0})= \left(
           \begin{array}{ccc}
             2      & 0         & 0 \\
             0      & 2         & 0 \\
             0      & 0         & 2 \\
           \end{array}
        \right),
\end{equation*}
and
\begin{equation}
    P(N_{q,0})=\frac{1}{2}A(N_{q,0}),
\end{equation}
respectively. The eigenvalues of $P(N_{q,0})$ are $1$ and $-\frac{1}{2}$, with their multiplicity being $1$ and $2$. Hence, by Lemmas~\ref{lemmaKem} and~\ref{lemmaK*}, the Kemeny constant and multiplicative degree-Kirchhoff index  for $N_{q,0}$  are $K(N_{q,0})=\frac{4}{3}$ and $\hat{\mathcal{K}}(N_{q,0})=8$, respectively.
Note that for the degree of each node in $N_{q,0}$ is 2.  By definition, for  $N_{q,0}$, its additive degree-Kirchhoff index is $\bar{\mathcal{K}}(N_{q,0})=8$, and its  Kirchhoff index is $\mathcal{K}(N_{q,0})=2$.
Then, by Theorems~\ref{iteratedKem},~\ref{iteratedK*},~\ref{iteratedK+}, and~\ref{iteratedKir}, we obtain the following exact solutions to the Kemeny's constant $K(N_{q,k})$, multiplicative degree-Kirchhoff index $\hat{\mathcal{K}}(N_{q,k})$, additive degree-Kirchhoff index $\bar{\mathcal{K}}(N_{q,k}) $, and Kirchhoff index $\mathcal{K}(N_{q,k})$ for  $N_{q,k}$.
\begin{small}
\begin{equation}\label{eg1}
\begin{aligned}
    K(N_{q,k})    =&    \frac{3(2q+3)(2q+1)^{k-1}}{2}
                        -\frac{q+4}{2q+1}\Big(\frac{4q+2}{q+2}\Big)^k\\
                  &     +\frac{4q+5}{6(2q+1)}.
\end{aligned}
\end{equation}
\end{small}
\begin{small}
\begin{equation}\label{eg2}
\begin{aligned}
  \hat{\mathcal{K}}(N_{q,k}) =&     \frac{9(2q+3)}{2q+1}(2q+1)^{2k}
                                    -\frac{6(q+4)}{2q+1}\Big(\frac{2(2q+1)^2}{q+2}\Big)^k\\
                              &     +(4q+5)(2q+1)^{k-1}.
\end{aligned}
\end{equation}
\end{small}
\begin{small}
\begin{equation}\label{eg4}
\begin{aligned}
        \bar{\mathcal{K}}(N_{q,k})
=&      \frac{9(2q+3)(6q+11)}{4(2q+1)(2q+5)}(2q+1)^{2k}\\
&       -\frac{3(q+4)}{2q+1}\Big(\frac{2(2q+1)^{2}}{q+2}\Big)^k\\
&       +\frac{3(q+4)}{2q+5}\Big(\frac{2(2q+1)}{q+2}\Big)^k\\
&       +\frac{4q+5}{2}(2q+1)^{k-1}
        +\frac{1}{4}.
\end{aligned}
\end{equation}
\end{small}
\begin{small}
\begin{equation}\label{eg5}
  \begin{aligned}
        \mathcal{K}(N_{q,k})
= &     \frac{9(2q+3)^2}{8(2q+1)(2q+5)}(2q+1)^{2k}\\
&       -\frac{3(q+4)}{8(2q+1)}\Big(\frac{2(2q+1)^{2}}{q+2}\Big)^k\\
&       +\frac{3(q+4)}{4(2q+5)}\Big(\frac{2(2q+1)}{q+2}\Big)^k\\
&       +\frac{(q+1)(4q+5)}{2(2q+5)}(2q+1)^{k-1}\\
&       +\frac{5(q+4)}{8(2q+5)}\Big(\frac{2}{q+2}\Big)^k
        -\frac{1}{8}.
  \end{aligned}
\end{equation}
\end{small}

\section{Conclusions}

The $q$-triangulation operation is an natural extension of traditional triangulation operation on a graph, which has been successfully applied to generate complex networks. In this paper, we presented an extensive study of various properties for $q$-triangulation graph $R_q(G)$ of a simple connected graph $G$, and obtained some interesting quantities of $R_q(G)$, which are expressed in terms of those associated with $G$. For this purpose, we first deduced formulas for eigenvalues and eigenvectors of normalized adjacency matrix of $R_q(G)$. Using these results, we then determined two-node hitting time and two-node resistance distance for an arbitrary node pair in $R_q(G)$. Also, we obtained the Kemeny's constant, Kirchhoff index, multiplicative degree-Kirchhoff index, and additive degree-Kirchhoff index for $R_q(G)$. As an application, we finally provided analytical formulas for some related quantities of iterated $q$-triangulations for a graph $G$, and obtained exact expressions for such quantities corresponding to pseudofractal scale-free webs, which mimic well realistic  networks with scale-free small-world properties.


%



\ifCLASSOPTIONcompsoc
  \section*{Acknowledgments}
\else
  \section*{Acknowledgment}
\fi

This work was supported by the National Natural Science Foundation of China under Grant No. 11275049.

\ifCLASSOPTIONcaptionsoff
  \newpage
\fi



%


\bibliographystyle{IEEEtran}
\bibliography{mybibfile,IEEEtran}

%

\begin{IEEEbiography}[{\includegraphics[width=1in,height=1.25in,clip,keepaspectratio]{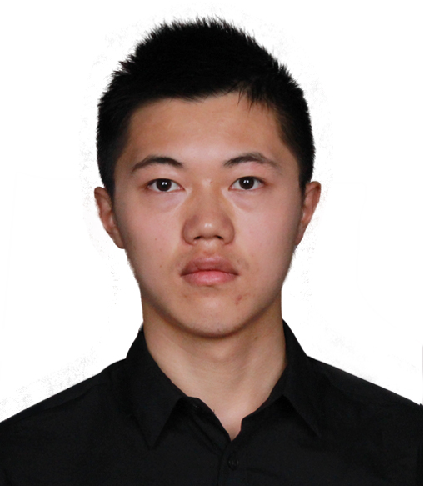}}]{Yibo Zeng}
is currently working toward the B.S. degree in the School of Mathematical Science, Fudan University. His research interests include network science, particularly in spectral properties of complex systems and networks.
\end{IEEEbiography}

\begin{IEEEbiography}[{\includegraphics[width=1in,height=1.25in,clip,keepaspectratio]{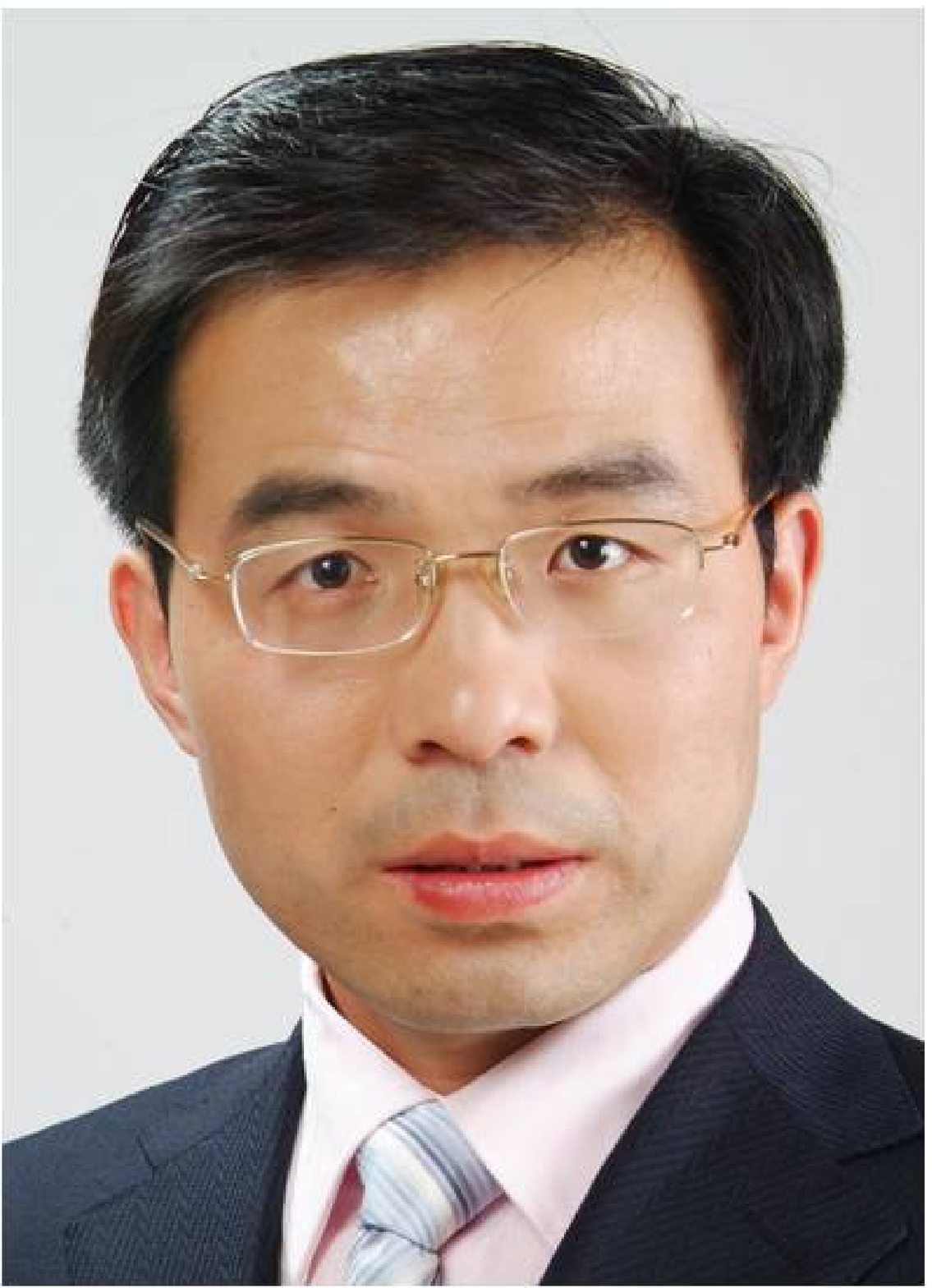}}]{Zhongzhi Zhang}
received the B.Sc. degree in Applied Mathematics from Anhui University, Hefei, China in 1997 and Ph.D. degree in Management Science and Engineering from Dalian University of Technology, Dalian, China, in 2006. From July 2006 to June 2008, he was a Post-Doctorate Research Fellow in Fudan University. Currently he is an Associate Professor at the School of Computer Science, Fudan University, China. Dr. Zhang's research interests include structural and dynamical properties of complex networks. He has published more than 100 papers in international journals since 2006 in the field of network modeling and dynamics, receiving over 1900 citations with H index 26 according to the Thomson Reuters ISI Web of Science.  Dr. Zhang is a Committee Member of the Chinese Society of Complexity Science. He received the Excellent Doctoral Dissertation Award of Liaoning Province, China in 2007, the Excellent Post-Doctor Award of Fudan University in 2008, and the Shanghai Natural Science Award (3rd class) in 2013.
\end{IEEEbiography}

%

%




\end{document}